\documentclass[12pt]{amsart}
\usepackage{amssymb}
\usepackage{mathrsfs}
\usepackage[all]{xy}

\RequirePackage{color}
\definecolor{myred}{rgb}{0.75,0,0}
\definecolor{mygreen}{rgb}{0,0.5,0}
\definecolor{myblue}{rgb}{0,0,0.65}

\RequirePackage{ifpdf}
\ifpdf
  \IfFileExists{pdfsync.sty}{\RequirePackage{pdfsync}}{}
  \RequirePackage[pdftex,
   colorlinks = true,
   urlcolor = myblue, 
   citecolor = mygreen, 
   linkcolor = myblue, 
   pdfusetitle,
   pdfauthor = {{Pramod Achar, Anthony Henderson, Daniel Juteau and Simon Riche}},
   pdfpagemode=None,
   bookmarksopen=true]{hyperref}
\else
  \RequirePackage[hypertex]{hyperref}
\fi

\usepackage{fullpage}

\usepackage{graphicx}
\usepackage{caption}
\usepackage{subcaption}
\usepackage{longtable}
\usepackage{supertabular}

\def\ov{\overline}

\newcommand{\C}{\mathbb{C}}
\newcommand{\Ql}{\mathbb{Q}_\ell}

\newcommand{\Qlb}{\overline{\mathbb{Q}}_\ell}
\newcommand{\Z}{\mathbb{Z}}

\newcommand{\K}{\mathbb{K}}

\renewcommand{\O}{\mathbb{O}}

\newcommand{\bk}{\Bbbk}

\newcommand{\fg}{\mathfrak{g}}

\newcommand{\cB}{\mathscr{B}}

\newcommand{\cN}{\mathscr{N}}
\newcommand{\cO}{\mathscr{O}}
\newcommand{\cC}{\mathscr{C}}
\newcommand{\cU}{\mathscr{U}}

\newcommand{\fN}{\mathfrak{N}}
\newcommand{\fM}{\mathfrak{M}}

\newcommand{\Db}{D^{\mathrm{b}}}
\newcommand{\Perv}{\mathsf{Perv}}

\newcommand{\Rep}{\mathsf{Rep}}

\newcommand{\For}{\mathrm{For}}
\newcommand{\bT}{\mathbb{T}}

\newcommand{\Res}{\mathbf{R}}
\newcommand{\Ind}{\mathbf{I}}
\newcommand{\IndRT}{\mathbf{Ind}}

\newcommand{\IC}{\mathcal{IC}}
\newcommand{\cF}{\mathcal{F}}

\newcommand{\ubk}{\underline{\bk}}
\newcommand{\Spr}{\underline{\mathsf{Spr}}}
\newcommand{\cE}{\mathcal{E}}

\newcommand{\simto}{\xrightarrow{\sim}}

\DeclareMathOperator{\Hom}{Hom}

\DeclareMathOperator{\Irr}{Irr}

\DeclareMathOperator{\Lie}{Lie}

\newcommand{\Ch}{\mathrm{Ch}}

\newcommand{\cusp}{\mathrm{cusp}}

\newcommand{\reg}{\mathrm{reg}}

\makeatletter
\def\lotimes{\@ifnextchar_{\@lotimessub}{\@lotimesnosub}}
\def\@lotimessub_#1{\mathchoice{\mathbin{\mathop{\otimes}^L}_{#1}}%
  {\otimes^L_{#1}}{\otimes^L_{#1}}{\otimes^L_{#1}}}
\def\@lotimesnosub{\mathbin{\mathop{\otimes}^L}}
\makeatother

\newcommand{\GL}{\mathrm{GL}}
\newcommand{\PGL}{\mathrm{PGL}}

\newcommand{\SL}{\mathrm{SL}}
\newcommand{\Sp}{\mathrm{Sp}}

\newcommand{\G}{{\mathrm{G}}}
\newcommand{\fS}{\mathfrak{S}}

\newcommand{\psico}{\psi^{\mathrm{co}}}

\newcommand{\omegaco}{\omega^{\mathrm{co}}}

\newcommand{\Part}{\mathrm{Part}}

\newcommand{\sm}{\mathsf{m}}
\newcommand{\uPart}{\underline{\Part}}

\newcommand{\blambda}{{\boldsymbol{\lambda}}}
\newcommand{\tr}{\mathrm{t}}

\newtheorem*{thm*}{Theorem}
\numberwithin{equation}{section}
\numberwithin{table}{section}
\newtheorem{thm}{Theorem}[section]

\newtheorem{prop}[thm]{Proposition}

\theoremstyle{definition}
\newtheorem{defn}[thm]{Definition}

\theoremstyle{remark}

\newtheorem{example}[thm]{Example}

\DeclareFontFamily{U}{mathx}{\hyphenchar\font45}
\DeclareFontShape{U}{mathx}{m}{n}{
      <5> <6> <7> <8> <9> <10>
      <10.95> <12> <14.4> <17.28> <20.74> <24.88>
      mathx10
      }{}
\DeclareSymbolFont{mathx}{U}{mathx}{m}{n}
\DeclareFontSubstitution{U}{mathx}{m}{n}
\DeclareMathAccent{\widebar}{0}{mathx}{"73}

\title[Modular generalized Springer correspondence]{Modular generalized Springer correspondence:\\ an overview}

\author{Pramod N. Achar}
\address{Department of Mathematics\\
  Louisiana State University\\
  Baton Rouge, LA 70803\\
  U.S.A.}
\email{pramod@math.lsu.edu}

\author{Anthony Henderson}
\address{School of Mathematics and Statistics\\
  University of Sydney, NSW 2006\\
  Australia}
\email{anthony.henderson@sydney.edu.au}

\author{Daniel Juteau}
\address{Laboratoire de Math\'ematiques Nicolas Oresme\\
  Universit\'e de Caen, BP 5186\\
  14032 Caen Cedex\\ 
  France}
\email{daniel.juteau@unicaen.fr}

\author{Simon Riche}
\address{Universit{\'e} Blaise Pascal - Clermont-Ferrand II, Laboratoire de Math{\'e}matiques, CNRS, UMR 6620, Campus universitaire des C{\'e}zeaux, F-63177 Aubi{\`e}re Cedex, France
}
\email{simon.riche@math.univ-bpclermont.fr}

\subjclass[2010]{Primary 17B08, 20G05}

\thanks{P.A. was supported by NSA Grant No.~H98230-15-1-0175 and NSF Grant No.~DMS-1500890.  A.H. was supported by ARC Future Fellowship Grant No.~FT110100504. D.J. and S.R. were supported by ANR Grant 
No.~ANR-13-BS01-0001-01. 
}

\begin{document}

\begin{abstract}
This is an overview of our series of papers on the modular generalized Springer correspondence~\cite{genspring1,genspring2,genspring3,rather-good}. It is an expansion of a lecture given by the second author in the Fifth Conference of the Tsinghua Sanya International Mathematics Forum, Sanya, December 2014, as part of the Master Lecture `Algebraic Groups and their Representations' Workshop honouring G.~Lusztig. The material that has not appeared in print before includes some discussion of the motivating idea of modular character sheaves, and heuristic remarks about geometric functors of parabolic induction and restriction.
\end{abstract}

\maketitle

\section{Motivation: modular character sheaves}
\label{sec:motivation}

One of Lusztig's fundamental contributions to the development of geometric representation theory was the theory of character sheaves on a connected reductive group~\cite{charsh1,charsh2,charsh3,charsh4,charsh5}. This theory complemented his earlier work with Deligne on the representation theory of finite reductive groups~\cite{delignelusztig,lusztig-book}, and was crucial to the later development of algorithms to compute the irreducible characters of such groups~\cite{char-val,green-fns,remarks,shoji-char-sheaves1,shoji-char-sheaves2}. Introductions to this theory are given in~\cite{charsh-intro,shoji-asterisque,marsspringer}.  

The subsequent thirty years have seen several generalizations of character sheaves, including generalizations to wider classes of algebraic groups; some of these are surveyed in~\cite{unity}, and see also~\cite{ginzburg,fgt,boyarchenko,boyarchenko-drinfeld,shoji-sorlin1,shoji-sorlin2,shoji-sorlin3}. To date, these generalizations, like Lusztig's original theory, have used sheaves or $D$-modules with coefficients in a field of characteristic zero, often $\Qlb$. Thus, they mimic ordinary (characteristic-zero) representations of finite groups. It is natural to hope for a theory of \emph{modular character sheaves} that mimics modular (positive characteristic) representations of finite groups. Such a theory should complement the existing modular applications of Deligne--Lusztig theory (see~\cite{gh} for a particularly relevant survey and~\cite{ce} for a general introduction and references), and could perhaps lead, for example, to a geometric interpretation of the decomposition matrix of a finite reductive group in the case of non-defining characteristic.

This hope builds on recent progress in other areas of \emph{modular geometric representation theory}, which may be loosely defined as the use of sheaves with characteristic-$\ell$ coefficients to answer questions about representations over fields of characteristic $\ell$. This theme arose in the work of Mirkovi\'c--Vilonen~\cite{mv} and Soergel~\cite{soergel}. Just as modular representation theory is made difficult by the failure of Maschke's Theorem, modular geometric representation theory is made difficult by the failure of the Decomposition Theorem of Be{\u\i}linson--Bernstein--Deligne--Gabber~\cite{bbd}. However, works such as~\cite{jmw} indicate ways to overcome this obstacle. For our present purposes, a key inspiration is the \emph{modular Springer correspondence} of the third author~\cite{juteau}, which affords a geometric interpretation of the decomposition matrix of a finite Weyl group by comparing $\Qlb$-sheaves and modular sheaves on the nilpotent cone.

We cannot yet formulate precise statements about, or even adequate definitions of, modular character sheaves. 
We will merely outline some features that such a theory may or may not be expected to possess, in comparison with Lusztig's theory.

In the following remarks, $G$ denotes a connected reductive group over an algebraically closed field $K$ of positive characteristic $p$, and $F:G\to G$ denotes a Frobenius endomorphism relative to some finite subfield of $K$, so that $G^F$ is a finite reductive group. We fix a prime $\ell$ different from $p$, and let $\bk$ be either $\Qlb$ or a finite extension of $\Z/\ell\Z$ (we refer to the latter possibility as the \emph{modular case}). For either choice of $\bk$, it makes sense to consider $\bk$-sheaves on $G$ for the \'etale topology.

\begin{enumerate}
\item \label{it:irreducibles}
In the $\bk=\Qlb$ case, Lusztig defined a collection $\widehat{G}$ of simple perverse $\Qlb$-sheaves on $G$, the \emph{character sheaves}. These correspond, in a somewhat loose sense recalled in~\eqref{it:Frobenius} below, to irreducible $\Qlb$-representations of the finite groups $G^F$ as $F$ varies. In the modular case, since we aim to mimic the non-semisimple categories $\Rep(G^F,\bk)$ of $\bk$-representations, we should define a category $\Ch(G,\bk)$ rather than just a collection of simple objects. To bring Lusztig's theory into the same framework, we could define $\Ch(G,\Qlb)$ to be the category of finite direct sums of elements of $\widehat{G}$.
\item \label{it:equivariance}
Lusztig, using the technology available at the time, defined character sheaves as certain objects $A$ of the bounded derived category $\Db(G,\Qlb)$, and more specifically of the subcategory $\Perv(G,\Qlb)$ of perverse sheaves, which are $G$-equivariant (for the $G$-action on itself by conjugation) in the sense that $a^*A\cong p^*A$ where $a:G\times G\to G$ is the action and $p:G\times G\to G$ is the second projection. This elementary notion of equivariance suffices for perverse sheaves, but the best definition of the equivariant derived category $\Db_G(G,\bk)$ is the one given subsequently by Bernstein and Lunts~\cite{bl}. It seems likely that $\Ch(G,\bk)$ should be defined as a 
subcategory of the abelian subcategory $\Perv_G(G,\bk)$ of $\Db_G(G,\bk)$, or at least (if some further structure is required) should have a forgetful functor to $\Perv_G(G,\bk)$.\footnote{
Another possibility might be to use parity sheaves as in~\cite{jmw}.
Following Soergel~\cite{soergel}, when the Decomposition Theorem fails we may consider the full additive subcategory of  $\Perv_G(G,\bk)$ generated by the direct summands of the direct images we are interested in, instead of the Serre 
subcategory generated by their subquotients. It would be nice to know that all those objects are parity sheaves, which according to \cite{jmw} can be classified in the same way as simple perverse sheaves. 
This would probably require proving first that the induction functors preserve parity complexes (which is known for the Springer morphism by \cite{dCLP} but not in the general case), and also the cleanness property in full generality (see Section \ref{sec:further}).
}
\item \label{it:induction}
As in the $\bk=\Qlb$ case, we would expect character sheaves to be stable under geometric functors of parabolic induction and restriction. We will illustrate what we mean in the case of induction. If $P$ is a parabolic subgroup of $G$ and $L$ is a Levi factor of $P$ so that $P=L\ltimes U_P$ where $U_P$ is the unipotent radical of $P$, then the parabolic induction functor
\[
\Ind_{L\subset P}^G:\Db_L(L,\bk)\to\Db_P(P,\bk)\to\Db_G(G,\bk)
\]
is the composition of a pull-back functor associated to the projection $P\twoheadrightarrow L$ and a push-forward functor associated to the inclusion $P\hookrightarrow G$.\footnote{See Section~\ref{sec:ind-res} below for a detailed definition of geometric parabolic induction in a slightly different setting.} Thus, in the case that $L$ and $P$ are $F$-stable, it is analogous to the representation-theoretic parabolic induction functor $\IndRT_{L^F\subset P^F}^{G^F}:\Rep(L^F,\bk)\to\Rep(G^F,\bk)$, also known as Harish--Chandra induction, which is the composition of the pull-back (inflation) functor associated to $P^F\twoheadrightarrow L^F$ and the induction functor associated to $P^F\hookrightarrow G^F$.
 We want  $\Ind_{L\subset P}^G : D^b_L(L,\bk)\to D^b_G(G,\bk)$ to induce a functor $\Ind_{L\subset P}^G : \Ch(L,\bk)\to\Ch(G,\bk)$, so that we can speak of cuspidal simple objects and induction series as usual (see Section~\ref{sec:cusp}). 
\item \label{it:Frobenius}
As the name `character sheaves' suggests, the relationship that Lusztig demonstrated between $\Ch(G,\Qlb)$ and $\Rep(G^F,\Qlb)$ goes via a matching-up of natural decategorifications of these categories, namely characteristic functions of sheaves and characters of representations. Recall that there are several subtleties in this matching-up. Firstly, one does not work with the whole of $\widehat{G}$ at once: rather, the choice of a Frobenius endomorphism $F$ determines a finite subset $\widehat{G}^F\subset\widehat{G}$ of $F$-stable character sheaves, those which are isomorphic to their pull-back under $F$. Secondly, the isomorphisms have to be suitably normalized in order to define for each $A\in\widehat{G}^F$ a characteristic function $\chi_A$ in $\cC(G^F,\Qlb)$, the vector space of $\Qlb$-valued class functions on $G^F$. Thirdly, $\{\chi_A\,|\,A\in\widehat{G}^F\}$ is not the set of characters of irreducible representations of $G^F$; rather, it turns out to be the set of \emph{almost} characters, a different basis of $\cC(G^F,\Qlb)$. 
Any relationship between $\Ch(G,\bk)$ and $\Rep(G^F,\bk)$ in the modular case will necessarily be even looser, because (as is well known) to get a reasonable decategorification of $\Rep(G^F,\bk)$ one must consider Brauer characters rather than the naive characters in $\cC(G^F,\bk)$.
It is not clear what will be the right way to take the Frobenius endomorphism into account in the modular case.
\item \label{it:analogy}
Despite the preceding remark, it might still be possible to draw analogies between operations on $\Ch(G,\bk)$ and operations on $\Rep(G^F,\bk)$, precisely by checking that they relate to the same operations on $\cC(G^F,\bk)$; in that sense, the name `character sheaves' is still justified in the modular case. For example, we can make precise the above suggestion that $\Ind_{L\subset P}^G:\Db_L(L,\bk)\to\Db_G(G,\bk)$ is analogous to $\IndRT_{L^F\subset P^F}^{G^F}:\Rep(L^F,\bk)\to\Rep(G^F,\bk)$. It is well known that the effect of the latter functor on the $\bk$-valued characters is given by the linear map $\mathrm{I}_{L^F\subset P^F}^{G^F}:\cC(L^F,\bk)\to\cC(G^F,\bk)$, defined by
\[
\mathrm{I}_{L^F\subset P^F}^{G^F}(f)(g)=\sum_{\substack{[g',p]\in G^F\times_{P^F} P^F\\g'p(g')^{-1}=g}}f(\overline{p}),
\]
where $\overline{\phantom{p}}:P^F\twoheadrightarrow L^F$ is the projection and $[g',p]$ denotes the orbit of a pair $(g',p)\in G^F\times P^F$ under the $P^F$-action defined by $p'\cdot(g',p)=(g'(p')^{-1},p'p(p')^{-1})$. It is indeed true in the modular case, just as in the $\bk=\Qlb$ case, that if $L$ and $P$ are $F$-stable, $A\in\Db_L(L,\bk)$ is $F$-stable, and the characteristic function of $A$ is $\chi_A\in\cC(G^F,\bk)$, then $\Ind_{L\subset P}^G(A)\in\Db_G(G,\bk)$ is $F$-stable and has characteristic function $\mathrm{I}_{L^F\subset P^F}^{G^F}(\chi_A)$. Although this statement may be less important in the modular setting (because $\bk$-valued characters are less useful than Brauer characters), it gives meaning to the analogy between $\Ind_{L\subset P}^G$ and $\IndRT_{L^F\subset P^F}^{G^F}$.
\end{enumerate}

Lusztig's definition of $\widehat{G}$ is purely geometric and applies equally well if the field $K$ over which $G$ is defined has characteristic zero. One could expect the same to be true of any theory of modular character sheaves. Of course, in this case $K$ has no finite subfields, so the analogy with representations of finite reductive groups is one step more remote. 

\section{The unipotent variety and the nilpotent cone}
\label{sec:unip-nilp}

Lusztig's theory of character sheaves is intimately related to his theory of the generalized Springer correspondence, contained in~\cite{lusztig,lusztig-fourier,lusztig-cusp1} and surveyed in~\cite{shoji-asterisque}. Here he studied $G$-equivariant perverse sheaves, not on $G$ itself, but on the unipotent variety $\cU_G\subset G$ and the nilpotent cone $\cN_G\subset\fg=\Lie(G)$. That is, he studied the abelian categories $\Perv_G(\cU_G,\Qlb)$ and $\Perv_G(\cN_G,\Qlb)$. Recall that $G$ has finitely many unipotent conjugacy classes and nilpotent adjoint orbits, so these categories have finitely many simple objects, unlike $\Perv_G(G,\Qlb)$; morally, this explains why there is no need to restrict attention to a proper subcategory like the subcategory $\Ch(G,\Qlb)\subset\Perv_G(G,\Qlb)$.

If the characteristic of the field of definition $K$ is either zero or sufficiently large relative to the root system of $G$, there is a $G$-equivariant isomorphism $\cN_G\simto\cU_G$ which induces an equivalence between $\Perv_G(\cU_G,\Qlb)$ and $\Perv_G(\cN_G,\Qlb)$. Thus, the difference between the two versions of the generalized Springer correspondence is negligible in a first approximation. 

Character sheaves relate to perverse sheaves on $\cU_G$ or $\cN_G$ in three interdependent ways, which we will do no more than mention here:
\begin{enumerate}
\item \label{it:restriction-G} \cite[Theorem 6.5]{lusztig} 
If $A\in\widehat{G}$ is a unipotent character sheaf,\footnote{By `unipotent character sheaf' we mean one of the simple perverse sheaves denoted $(\phi_!K)_\rho$ in~\cite[Theorem 6.5(c)]{lusztig}, which were shown to be character sheaves in~\cite{charsh5}.} the restriction $A|_{\cU_G}$ is a shift of a simple object of $\Perv_G(\cU_G,\Qlb)$.
\item \label{it:restriction-g} \cite{lusztig-fourier} 
There is a collection of character sheaves\footnote{In Lusztig's terminology in~\cite{lusztig-fourier}, these are the `admissible' objects in $\Perv_G(\fg,\Qlb)$. The statements recalled here are special cases of general results about admissible objects,~\cite[3(b) and Theorem 5(a)]{lusztig-fourier}.} on $\fg$ analogous to the character sheaves on $G$. If $A$ is a nilpotent character sheaf on $\fg$ (the analogue of a unipotent character sheaf on $G$), then not only is the restriction $A|_{\cN_G}$ a shift of a simple object of $\Perv_G(\cN_G,\Qlb)$, the Fourier--Deligne transform $\bT_\fg(A)$ of $A$ is itself a (generally different) simple object of $\Perv_G(\cN_G,\Qlb)$, extended by zero to the whole of $\fg$.
\item \label{it:extension-by-zero} \cite{lusztig,lusztig-fourier} 
Assume for convenience in this statement that $G$ is semisimple. Then cuspidal simple objects of $\Perv_G(\cU_G,\Qlb)$ (respectively, $\Perv_G(\cN_G,\Qlb)$), when extended by zero to the whole of $G$ (respectively, $\fg$), are important examples of cuspidal character sheaves.
\end{enumerate} 

For all these reasons, in working towards a theory of modular character sheaves, a natural first step is a modular version of the generalized Springer correspondence, in which one considers $\Perv_G(\cU_G,\bk)$ and $\Perv_G(\cN_G,\bk)$ for the more general fields $\bk$ introduced in the previous section. 

In the modular case, there is a strong reason to prefer $\cN_G\subset\fg$ to $\cU_G\subset G$. Namely,
there seems to be little hope of a modular analogue of the restriction statements in \eqref{it:restriction-G} and \eqref{it:restriction-g}, since 
the constant $\bk$-sheaf $\ubk$ on $\cU_G$ or $\cN_G$ is not a shift of a simple perverse sheaf in general (see e.g.~\cite[Section 8]{juteau}). However, the Fourier--Deligne statement in \eqref{it:restriction-g} 
seems more likely to have a modular analogue.
This phenomenon was encountered by the third author in the development of the modular Springer correspondence~\cite{juteau}, in which Fourier--Deligne transform is the vital tool. 

Henceforth, we restrict attention to $\Perv_G(\cN_G,\bk)$. Moreover, we make the simplifying assumption that $K=\C$; this allows us to isolate the new phenomena associated with the change of coefficients from the $\bk=\Qlb$ case to the modular case, without worrying about the phenomena associated with small values of $p$.\footnote{
For example, when $p$ is a bad prime for $G$, the classification of nilpotent orbits and associated structures are different; some of our proofs use case-by-case considerations which would not automatically carry over to this situation.
}

Largely for technical convenience, we consider $\bk$-sheaves on $\cN_G$ for the complex topology rather than the \'etale topology.\footnote{
This is mostly to avoid introducing the \'etale and $\ell$-adic constructible equivariant derived categories; once that formalism is in place, exactly the same arguments work.}
We can then allow $\bk$ to be any field (or even, for some purposes, a commutative ring as in~\cite{small2,ahjr-weyl}). The replacement for Fourier--Deligne transform in this setting is the Fourier--Sato transform defined in~\cite{ks}.

\section{Parabolic induction and restriction}
\label{sec:ind-res}

We have now arrived at the setting of our series of papers~\cite{genspring1,genspring2,genspring3,rather-good}. For clarity, we restate the assumptions that apply henceforth:
\begin{itemize}
\item $G$ is a connected reductive group over $\C$;
\item $\cN_G$ is the nilpotent cone in $\fg=\Lie(G)$;
\item $\bk$ is a field of characteristic $\ell\geq 0$ (where the $\ell=0$ case can be regarded as already known from Lusztig's work);
\item $\Db_G(\cN_G,\bk)$, $\Perv_G(\cN_G,\bk)$ are defined as in~\cite{bl} using the complex topology on $\cN_G$. 
\end{itemize}
We furthermore assume that $\bk$ is big enough for the following condition to hold:
\begin{equation} \label{eqn:big-enough}
\begin{array}{c}
\text{$\bk$ is a splitting field for all the finite groups $L_x/L_x^\circ$}\\
\text{where $L$ is a Levi subgroup of $G$ and $x\in\cN_L$.}
\end{array}
\end{equation} 
Here $L_x$ is the stabilizer of $x$ in $L$, and $L_x/L_x^\circ$ is its group of components.
If $\bk$ is algebraically closed,~\eqref{eqn:big-enough} certainly holds; see~\cite[Section 3.2]{genspring3} for a detailed analysis of~\eqref{eqn:big-enough} for various groups $G$, which shows that it is a very mild assumption.

Let $P$ be a parabolic subgroup of $G$ and $L$ a Levi factor of $P$. Our first task is to define a parabolic induction functor $\Ind_{L\subset P}^G:\Db_L(\cN_L,\bk)\to\Db_G(\cN_G,\bk)$.  
As was motivated above in the context of character sheaves, this should be the composition of a pull-back functor associated to the projection $P\twoheadrightarrow L$ and a push-forward functor associated to the inclusion $P\hookrightarrow G$, but we need to interpret these appropriately for equivariant derived categories. We will explain our thinking in some detail, to provide a model for possible generalization.
 
Recall from~\cite[Section 6]{bl} that one has an equivariant pull-back functor $Q_f^*$ and a right adjoint equivariant push-forward functor $Q_{f*}$ for any map $f:X\to Y$, where the groups acting on $X$ and $Y$ are $G$ and $H$ respectively, and $f$ is $G$-equivariant for the $G$-action on $Y$ defined using a given homomorphism $\phi:G\to H$. The push-forward functor in general only preserves the equivariant \emph{bounded-below} derived categories, but we will consider certain special cases where it does restrict to a functor $Q_{f*}:\Db_G(X,\bk)\to\Db_H(Y,\bk)$. The pull-back functor $Q_f^*:\Db_H(Y,\bk)\to\Db_G(X,\bk)$ always preserves the equivariant bounded derived categories. These categories possess a Verdier duality $D$ by~\cite[Theorem 3.5.2]{bl}, so we can also set $Q_f^!=D\circ Q_f^*\circ D$ and $Q_{f!}=D\circ Q_{f*}\circ D$, with $Q_f^!$ being right adjoint to $Q_{f!}$.

From the stack viewpoint, $Q_f^*$, $Q_{f*}$, $Q_f^!$, $Q_{f!}$ are just the standard four functors for the induced map of quotient stacks $[G\backslash X]\to[H\backslash Y]$. Since the latter map factors as $[G\backslash X]\to[G\backslash Y]\to[H\backslash Y]$, we can express each of the four functors as a composition where one factor relates $\Db_G(X,\bk)$ and $\Db_G(Y,\bk)$ and the other factor relates $\Db_G(Y,\bk)$ and $\Db_H(Y,\bk)$. Henceforth we will do so, at the risk of obscuring the essential simplicity of the definitions. 

Applying this formalism to the inclusion $i_{L\subset P}:\cN_P\hookrightarrow\cN_G$ of varieties and the associated inclusion $P\hookrightarrow G$ of groups, we obtain the pull-back functor $(i_{L\subset P})^*\circ\For_P^G$, where $(i_{L\subset P})^*:\Db_P(\cN_G,\bk)\to\Db_P(\cN_P,\bk)$ is as defined in~\cite[Section 3]{bl} and $\For_P^G:\Db_G(\cN_G,\bk)\to\Db_P(\cN_G,\bk)$ is the forgetful functor defined in~\cite[Section 2.6.1]{bl} (where it was denoted $Res_{P,G}$). The right adjoint push-forward functor is $\Gamma_P^G\circ (i_{L\subset P})_*$, where $(i_{L\subset P})_*:\Db_P(\cN_P,\bk)\to\Db_P(\cN_G,\bk)$ is as defined in~\cite[Section 3]{bl} and $\Gamma_P^G:\Db_P(\cN_G,\bk)\to\Db_G(\cN_G,\bk)$ is the right adjoint to $\For_P^G$ defined in~\cite[Theorem 3.7.1]{bl}. The Verdier-dual versions of these functors are $(i_{L\subset P})^!\circ\For_P^G$ and $\gamma_P^G\circ (i_{L\subset P})_!$, where $\gamma_P^G$ is the left adjoint to $\For_P^G$, discussed in~\cite[Section B.10.1]{small2}. Since $i_{L\subset P}$ is a closed embedding and $G/P$ is complete,
\[
\gamma_P^G\circ (i_{L\subset P})_!\cong\Gamma_P^G\circ (i_{L\subset P})_*[2\dim(G/P)],
\]
so a Verdier-self-dual version of the push-forward functor is $\gamma_P^G\circ (i_{L\subset P})_![-\dim(G/P)]$.

Applying the same formalism to the projection $p_{L\subset P}:\cN_P\twoheadrightarrow\cN_L$ and the associated projection $P\twoheadrightarrow L$, we obtain the pull-back functor $(p_{L\subset P})^*\circ \gamma_L^P$, because $\gamma_L^P\cong\Gamma_L^P:\Db_L(\cN_L,\bk)\to\Db_P(\cN_L,\bk)$ and $\For_L^P:\Db_P(\cN_L,\bk)\to\Db_L(\cN_L,\bk)$ are inverse equivalences by~\cite[Theorem 3.7.3]{bl}. The right adjoint push-forward functor is $\For_L^P\circ (p_{L\subset P})_*$. The Verdier-dual versions of these functors are $(p_{L\subset P})^!\circ \Gamma_L^P$ and $\For_L^P\circ (p_{L\subset P})_!$. Since $p_{L\subset P}$ is smooth with fibres of dimension $\dim(G/P)$, we have
\[
(p_{L\subset P})^*\circ \gamma_L^P\cong(p_{L\subset P})^!\circ \Gamma_L^P[-2\dim(G/P)],
\]
so a Verdier-self-dual version of the pull-back functor is $(p_{L\subset P})^*\circ \gamma_L^P[\dim(G/P)]$.    

The upshot of these comments is the following definition of parabolic induction:
\[
\Ind_{L\subset P}^G=\gamma_P^G\circ(i_{L\subset P})_!\circ(p_{L\subset P})^*\circ \gamma_L^P\cong\Gamma_P^G\circ(i_{L\subset P})_*\circ(p_{L\subset P})^!\circ \Gamma_L^P:\Db_L(\cN_L,\bk)\to\Db_G(\cN_G,\bk).
\]
We stated the first of these equivalent definitions in~\cite[Section 2.1]{genspring1}.\footnote{In~\cite[Section 2.1]{genspring1} we omitted $\gamma_L^P$, since it can be thought of as an identification.} The fact that the two definitions are equivalent says that $\Ind_{L\subset P}^G$ commutes with Verdier duality. See~\cite[Lemma 2.14]{genspring1} for yet another equivalent definition of $\Ind_{L\subset P}^G$, which is closer to the form that Lusztig used.

We see immediately that $\Ind_{L\subset P}^G$ has right and left adjoint functors of parabolic restriction, respectively
\[
\begin{split}
\Res_{L\subset P}^G&=\For_L^P\circ(p_{L\subset P})_*\circ(i_{L\subset P})^!\circ\For_P^G:\Db_G(\cN_G,\bk)\to\Db_L(\cN_L,\bk),\\
{}'\Res_{L\subset P}^G&=\For_L^P\circ(p_{L\subset P})_!\circ(i_{L\subset P})^*\circ\For_P^G:\Db_G(\cN_G,\bk)\to\Db_L(\cN_L,\bk),
\end{split}
\]
which are interchanged by Verdier duality.\footnote{In stating the definitions of $\Res_{L\subset P}^G$ and ${}'\Res_{L\subset P}^G$ in~\cite[Section 2.1]{genspring1}, we suppressed the forgetful functors.}

The following key properties of these induction and restriction functors, needed to get the theory started, are proved in~\cite[Proposition 4.7, Corollary 4.13]{small2} and~\cite[Proposition 3.1]{am}, using Braden's hyperbolic restriction theorem~\cite{braden}.

\begin{prop} \label{prop:adjunctions}
$\Ind_{L\subset P}^G$, $\Res_{L\subset P}^G$ and ${}'\Res_{L\subset P}^G$ restrict to exact functors between the perverse subcategories $\Perv_L(\cN_L,\bk)$ and $\Perv_G(\cN_G,\bk)$, which of course still satisfy the adjunctions
\[
{}'\Res_{L\subset P}^G \vdash \Ind_{L\subset P}^G \vdash \Res_{L\subset P}^G.
\] 
Moreover, we have an isomorphism
\[
{}'\Res_{L\subset P}^G\cong\Res_{L\subset P^-}^G\text{ as functors }\Perv_G(\cN_G,\bk)\to\Perv_L(\cN_L,\bk),
\] 
where $P^-$ denotes the parabolic subgroup opposite to $P$ with the same Levi factor $L$.
\end{prop}

It may seem surprising that we do not have an isomorphism ${}'\Res_{L\subset P}^G\cong\Res_{L\subset P}^G$ and hence a biadjunction between $\Ind_{L\subset P}^G$ and $\Res_{L\subset P}^G$. In Lusztig's setting, one can use the fact that $\Perv_G(\cN_G,\Qlb)$ and $\Perv_L(\cN_L,\Qlb)$ are semisimple to prove such a biadjunction; this proof does not apply in the modular setting, since $\Perv_G(\cN_G,\bk)$ is semisimple only when $\ell$ does not divide the order of the Weyl group $W$ of $G$, see~\cite[Proposition 7.9]{rather-good}. 

In the representation theory of finite reductive groups (with $\ell\neq p$ as in Section~\ref{sec:motivation}), one does have a biadjunction between parabolic induction and restriction. In a sense, this relies on the fact that in characteristic $\ell$, invariants and coinvariants for the finite $p$-group $U_P^F$ are the same. There is no obvious geometric analogue of this fact.

\section{Cuspidal pairs and induction series}
\label{sec:cusp}

Let $\fN_{G,\bk}$ denote the set of pairs $(\cO,\cE)$ where $\cO$ is a $G$-orbit in $\cN_G$ and $\cE$ runs over the irreducible $G$-equivariant $\bk$-local systems on $\cO$ (taken up to isomorphism). Recall that such irreducible local systems $\cE$ on $\cO$ are in a natural bijection with irreducible $\bk$-representations of the finite group $A_G(x)=G_x/G_x^\circ$, 
where $x\in\cO$ is a chosen point.
The assumption~\eqref{eqn:big-enough} ensures that each such $\cE$ is absolutely irreducible, and that the same holds for all Levi subgroups of $G$.

The isomorphism classes of simple objects in $\Perv_G(\cN_G,\bk)$ are in bijection with $\fN_{G,\bk}$. For $(\cO,\cE)\in\fN_{G,\bk}$, the corresponding simple perverse sheaf is the intersection cohomology complex $\IC(\cO,\cE)$, also known as the intermediate extension ${}^p j_{!*}(\cE[\dim\cO])$ where $j:\cO\hookrightarrow\overline{\cO}$ is the inclusion, extended by zero from $\overline{\cO}$ to $\cN$. See~\cite{bbd} for the definition of $\IC(\cO,\cE)$, and~\cite{jmw1} for some worked examples showing the difference between the characteristic-$0$ and modular cases.   

\begin{defn} \label{def:cuspidal} \cite[Definition 2.2]{genspring1}
We say that a pair $(\cO,\cE)\in\fN_{G,\bk}$, or the corresponding simple perverse sheaf $\IC(\cO,\cE)$, is \emph{cuspidal} if the following equivalent conditions hold (where $P$ and $L$ denote parabolic and Levi subgroups as before):
\begin{enumerate}
\item $\Res_{L\subset P}^G(\IC(\cO,\cE))=0$ for all $L\subset P\subsetneq G$;
\item ${}'\Res_{L\subset P}^G(\IC(\cO,\cE))=0$ for all $L\subset P\subsetneq G$;
\item $\IC(\cO,\cE)$ is not a quotient of $\Ind_{L\subset P}^G(A)$ for any $L\subset P\subsetneq G$ and any $A\in\Perv_L(\cN_L,\bk)$;
\item $\IC(\cO,\cE)$ is not a subobject of $\Ind_{L\subset P}^G(A)$ for any $L\subset P\subsetneq G$ and any $A\in\Perv_L(\cN_L,\bk)$.
\end{enumerate}
The equivalence of these conditions follows immediately from Proposition~\ref{prop:adjunctions}, see~\cite[Proposition~2.1]{genspring1}. We let $\fN_{G,\bk}^\cusp\subset\fN_{G,\bk}$ be the set of cuspidal pairs.
\end{defn}

When $\ell=0$, the Decomposition Theorem of~\cite{bbd} implies that if $A\in\Perv_L(\cN_L,k)$ is simple, then $\Ind_{L\subset P}^G(A)$ is semisimple, so one can replace `quotient' or `subobject' in Definition~\ref{def:cuspidal} by `simple constituent'. When $\ell>0$, the failure of the Decomposition Theorem means that $\Ind_{L\subset P}^G(A)$ is sometimes not semisimple, and it can have simple constituents (even cuspidal ones) which are not quotients or subobjects; see~\cite[Section 8]{juteau} for an example in the case $G=\mathrm{SL}(2)$ and $\ell=2$. Of course, this situation is familiar in modular representation theory.

A key ingredient in Lusztig's treatment of cuspidal pairs is that in the $\ell=0$ case there is another criterion for cuspidality, equivalent to those in Definition~\ref{def:cuspidal}, which refers only to the local system $\cE$ and not to the perverse sheaf $\IC(\cO,\cE)$~\cite[Proposition 2.2]{lusztig}. We showed in~\cite[Proposition 2.4]{genspring1} that in the $\ell>0$ case Lusztig's criterion implies cuspidality but is not implied by it. In view of this, one can immediately expect that there will be more cuspidal pairs in the modular case, and that they will be harder to classify. However, one of Lusztig's results about cuspidal pairs does generalize easily:

\begin{prop} \label{prop:distinguished} \cite[Proposition 2.6]{genspring2}
If $(\cO,\cE)$ is cuspidal, then $\cO$ is distinguished. 
\end{prop}
\noindent
Recall that a nilpotent orbit is said to be \emph{distinguished} if it does not meet the Lie algebra of any proper Levi subgroup $L\subsetneq G$.

A \emph{cuspidal datum} for $G$ is a triple $(L,\cO_L,\cE_L)$ where $L$ is a Levi subgroup of $G$ (possibly equal to $G$) and $(\cO_L,\cE_L)\in\fN_{L,\bk}^\cusp$ is a cuspidal pair for $L$. We let $\fM_{G,\bk}$ denote either the set of $G$-conjugacy classes of cuspidal data or, in a slight abuse of notation, a set of representatives of these $G$-conjugacy classes. 

Any $(L,\cO_L,\cE_L)\in\fM_{G,\bk}$ gives rise to a parabolically induced perverse sheaf
\[
\Ind(L,\cO_L,\cE_L)=\Ind_{L\subset P}^G(\IC(\cO_L,\cE_L))\in\Perv_G(\cN_G,\bk),
\]
where $P$ is a parabolic subgroup of $G$ having $L$ as a Levi factor.

\begin{prop} \label{prop:independence} \cite[Section 2]{genspring2}
Up to isomorphism, $\Ind(L,\cO_L,\cE_L)$ is independent of the choice of parabolic subgroup $P$, and depends only on the $G$-conjugacy class of $(L,\cO_L,\cE_L)$. The head and socle of $\Ind(L,\cO_L,\cE_L)$ are isomorphic to each other.
\end{prop}

Proposition~\ref{prop:independence} is an analogue of a known result in modular representation theory~\cite[Section 2.2]{gh}. However, unlike in that context, we do not prove that $\Ind_{L\subset P}^G(A)$ is independent of $P$ for all $A\in\Perv_L(\cN_L,\bk)$, only for cuspidal simple $A$. Our proof proceeds by describing the Fourier--Sato transform of $\Ind(L,\cO_L,\cE_L)$ as an intersection cohomology complex on $\fg$: as in the $\ell=0$ case, it is the intermediate extension of a certain local system on a piece of the Lusztig stratification of $\fg$, defined in~\cite[Section 6]{lusztig-cusp2}.  

\begin{defn} \label{def:induction-series}
For $(L,\cO_L,\cE_L)\in\fM_{G,\bk}$, the \emph{induction series} associated to $(L,\cO_L,\cE_L)$ is the set of isomorphism classes of simple quotients (equivalently, simple subobjects) of $\Ind(L,\cO_L,\cE_L)$, or the corresponding set of pairs $\fN_{G,\bk}^{(L,\cO_L,\cE_L)}\subset\fN_{G,\bk}$.
\end{defn}

As in the $\ell=0$ case, it is a straightforward consequence of the transitivity of parabolic induction~\cite[Lemma 2.6]{genspring1} that 
\[
\fN_{G,\bk} = \bigcup_{(L,\cO_L,\cE_L) \in \fM_{G,\bk}} \fN_{G,\bk}^{(L,\cO_L,\cE_L)},
\]
see~\cite[Corollary~2.7]{genspring1}.
In more concrete terms, every simple perverse sheaf $\IC(\cO,\cE)$ in $\Perv_G(\cN_G,\bk)$ occurs as a quotient of $\Ind(L,\cO_L,\cE_L)$ for some $(L,\cO_L,\cE_L)\in\fM_{G,\bk}$.

At one extreme, we have cuspidal data of the form $(G,\cO,\cE)$ where $(\cO,\cE)\in\fN_{G,\bk}^\cusp$: the associated induction series consists just of $(\cO,\cE)$. At the other extreme, we have the cuspidal datum $(T,\{0\},\ubk)$ where $T$ is a maximal torus of $G$: the associated induction series is called the \emph{principal series}.  By~\cite[Lemma 2.14]{genspring1}, the parabolically induced perverse sheaf $\Ind(T,\{0\},\ubk)$ is exactly the \emph{Springer sheaf} $\Spr=\mu_*\ubk[\dim\cN_G]$, where $\mu:T^*\cB\to\cN_G$ is Springer's resolution of the nilpotent cone, so the principal series can alternatively be defined as the set of isomorphism classes of simple quotients of the Springer sheaf. 

The \emph{modular Springer correspondence} of the third author~\cite{juteau} is a canonical bijection between the principal series $\fN_{G,\bk}^{(T,\{0\},\ubk)}$ and the set $\Irr(\bk[W])$ of isomorphism classes of irreducible $\bk$-representations of the Weyl group $W$ of $G$, generalizing Springer's original correspondence in the $\ell=0$ case. The group $W$ has a natural action\footnote{Actually there are two ways to define the $W$-action, using either restriction or Fourier--Sato transform as in Section~\ref{sec:unip-nilp}, but it is proved in~\cite{ahjr-weyl} that they differ only by tensoring with the sign representation.} on $\Spr$ by automorphisms in $\Perv_G(\cN_G,\bk)$, and the modular Springer correspondence sends a pair $(\cO,\cE)\in\fN_{G,\bk}^{(T,\{0\},\ubk)}$ to the representation $\Hom(\Spr,\IC(\cO,\cE))$. This correspondence has been explicitly described in the various cases in~\cite{juteau,jls,genspring2}.

\section{The main theorem}
\label{sec:mgsc}

We are now ready to state our main theorem, the \emph{modular generalized Springer correspondence}, which generalizes the aforementioned modular Springer correspondence, and includes Lusztig's generalized Springer correspondence as the $\ell=0$ case. The general statement is~\cite[Theorem 1.1]{genspring3}; special cases appeared separately as~\cite[Theorem 3.3]{genspring1} (the case where $G=\GL(n)$) and~\cite[Theorem 1.1]{genspring2} (the case where $G$ is a classical group).

\begin{thm} \label{thm:mgsc}
We keep the assumptions made at the start of Section~{\rm \ref{sec:ind-res}}. 
\begin{enumerate}
\item
\label{it:mgsc-1}
The induction series associated to different cuspidal data are disjoint, so that we have
\[
\fN_{G,\bk} = \bigsqcup_{(L,\cO_L,\cE_L) \in \fM_{G,\bk}} \fN_{G,\bk}^{(L,\cO_L,\cE_L)}.
\]
In other words, every simple perverse sheaf $\IC(\cO,\cE)\in\fN_{G,\bk}$ occurs as a quotient of $\Ind(L,\cO_L,\cE_L)$ for a \emph{unique} $(L,\cO_L,\cE_L)\in\fM_{G,\bk}$.
\item
\label{it:mgsc-2}
For any $(L,\cO_L,\cE_L)\in\fM_{G,\bk}$, we have a canonical bijection between the induction series $\fN_{G,\bk}^{(L,\cO_L,\cE_L)}$ and the set $\Irr(\bk[N_G(L)/L])$ of isomorphism classes of irreducible $\bk$-representations of $N_G(L)/L$.
\item
\label{it:mgsc-3}
Hence, combining~\eqref{it:mgsc-1} and~\eqref{it:mgsc-2}, we have a bijection
\[ \fN_{G,\bk}\;\longleftrightarrow \bigsqcup_{(L,\cO_L,\cE_L)\in\fM_{G,\bk}} \Irr(\bk[N_G(L)/L]), \]
called the modular generalized Springer correspondence.
\end{enumerate}
\end{thm}

The bijection defined in~\eqref{it:mgsc-2} is independent of the chosen representative $(L,\cO_L,\cE_L)$ of a $G$-conjugacy class in $\fM_{G,\bk}$, in the following sense. If $(M,\cO_M,\cE_M) = g \cdot (L,\cO_L,\cE_L)$ for some $g \in G$, then conjugation by $g$ induces an isomorphism $N_G(L)/L \simto N_G(M) / M$. The induced bijection $\Irr(\bk[N_G(L)/L]) \leftrightarrow \Irr(\bk[N_G(M)/M])$ does not depend on the choice of $g$, and its composition with the bijection defined in~\eqref{it:mgsc-2} is the analogous bijection for $(M,\cO_M,\cE_M)$.

The disjointness of induction series is well known in the context of the representation theory of finite reductive groups, and, as in that theory, our general proof of part~\eqref{it:mgsc-1} uses a Mackey formula for the composition of functors $\Res_{M\subset Q}^G\circ\Ind_{L\subset P}^G$. In our setting, this Mackey formula~\cite[Theorem 2.2]{genspring3} has a filtration rather than a direct sum. It is analogous to the Mackey formula for character sheaves proved in~\cite[Section 10.1]{marsspringer}, and the Mackey formula for $D$-modules on Lie algebras proved in~\cite[Section 3]{gunningham}. Our proof in~\cite{genspring3} imitates the proof of the former, but could have been phrased in the language of stacks like the proof of the latter.   

Part~\eqref{it:mgsc-2} differs from the analogous statement in representation theory in two ways; these differences were already apparent in Lusztig's result in the $\ell=0$ case. First, we have the full group $N_G(L)/L$ rather than a subgroup, because the normalizer $N_G(L)$ preserves $(\cO_L,\cE_L)$; this can be deduced, using Proposition~\ref{prop:distinguished}, from a general result on distinguished orbits~\cite[Proposition~3.1]{genspring3}. Second, we have the group algebra $\bk[N_G(L)/L]$ rather than a deformation such as a Hecke algebra; this is a more delicate matter.
Our proof of this aspect, which is contained in~\cite[Section 3]{genspring2}, uses the description of the Fourier--Sato transform of $\Ind(L,\cO_L,\cE_L)$, mentioned above in connection with Proposition~\ref{prop:independence}. It relies on a detailed analysis of certain resolutions of pieces of the Lusztig stratification of $\fg$ defined by Bonnaf\'e~\cite{bonnafe1}; it is this geometry which allows us to show that a certain endomorphism algebra, which a priori is a \emph{twisted} group algebra of $N_G(L)/L$, is in fact \emph{canonically} isomorphic to $\bk[N_G(L)/L]$.

In this proof we are required to diverge from Lusztig's proofs in~\cite{lusztig} and~\cite{lusztig-cusp2}, mainly because of one key difference from his setting. In the $\ell=0$ case, Lusztig shows in~\cite[Theorem 9.2]{lusztig} that any Levi subgroup $L$ supporting a cuspidal pair must be \emph{self-opposed} in $G$ in the terminology of~\cite[Section 1.E]{bonnafe1}, a serious restriction which implies in particular that $N_G(L)/L$ is a finite Weyl group. The trivial and sign representations of $N_G(L)/L$ then play important special roles in the generalized Springer correspondence.\footnote{\label{footnote:sign-twist}The bijection between $\fN_{G,\bk}^{(L,\cO_L,\cE_L)}$ and $\Irr(\bk[N_G(L)/L])$ defined by Lusztig in his setting~\cite[Theorem 6.5(c) and Theorem 9.2(d)]{lusztig} differs from the $\ell=0$ case of our bijection by tensoring with the sign representation of $N_G(L)/L$. This is because he relates character sheaves on $G$ to perverse sheaves on $\cU_G$ by restriction (see Section~\ref{sec:unip-nilp} above), whereas we work with $\cN_G\subset\fg$ and use Fourier--Sato transform.} This is a key ingredient not only of Lusztig's proof of part~\eqref{it:mgsc-2}, but also of Spaltenstein's explicit determination of the generalized Springer correspondence for exceptional groups~\cite{spaltenstein}. 

In the modular case, a Levi subgroup $L$ supporting a cuspidal pair need not be self-opposed in $G$, as one sees already in the case where $G=\GL(n)$ (see Example~\ref{ex:gln} below). In fact, the group $N_G(L)/L$ need not even be a reflection group~\cite[Remark 6.5]{genspring3}. The situation appears to be much less rigid than the $\ell=0$ case: for instance, we have not found a general rule predicting which element of the induction series $\fN_{G,\bk}^{(L,\cO_L,\cE_L)}$ corresponds to the trivial representation of $N_G(L)/L$, analogous to Lusztig's result~\cite[Proposition 9.5]{lusztig}.\footnote{Lusztig's result~\cite[Proposition 9.5]{lusztig} concerns the sign representation, but that becomes the trivial representation in our conventions, as per Footnote~\ref{footnote:sign-twist}. For example, in our conventions, the pair in the principal series $\fN_{G,\bk}^{(T,\{0\},\ubk)}$ corresponding to the trivial representation of $W$ is $(\{0\},\ubk)$.}

Finally, let us mention that one can in some cases obtain statements~\eqref{it:mgsc-1} and~\eqref{it:mgsc-2} by simpler methods. For instance, \eqref{it:mgsc-2} is easy if every $L$-equivariant local system on $\cO_L$ is constant, and~\eqref{it:mgsc-1} is easy if one knows that no proper Levi subgroup $L$ of $G$ has two cuspidal pairs supported on the same orbit with the same central character (see~\cite[Corollary~2.3]{genspring2}). In particular these conditions clearly hold when $G=\mathrm{GL}(n)$ (because in this case all equivariant local systems on all nilpotent orbits of all Levi subgroups are constant -- see Example~\ref{ex:gln} below), which allowed us to treat this case in~\cite{genspring1}, before developing the general theory. More generally, for classical groups, we were able to check the second condition by explicitly classifying the cuspidal pairs (at the same time as establishing the correspondence) in~\cite{genspring2}; hence in that paper we did not need the Mackey formula.

\section{Examples}
\label{sec:examples}

\begin{example} \label{ex:gln}
Suppose that $G=\GL(n)$ for some positive integer $n$; this case is treated in detail in~\cite{genspring1}. The $G$-orbits in $\cN_G$ are in bijection with the set $\Part(n)$ of partitions of $n$: for a partition $\lambda=(\lambda_1,\lambda_2,\cdots)\in\Part(n)$, the corresponding orbit $\cO_\lambda$ consists of nilpotent matrices whose Jordan blocks have sizes $\lambda_1,\lambda_2,\cdots$. The unique distinguished nilpotent orbit is the regular orbit $\cO_{\reg}=\cO_{(n)}$.

In this case $A_G(x)=1$ for all $x\in\cN_G$, so all pairs in $\fN_{G,\bk}$ have the form $(\cO_\lambda,\ubk)$ for some $\lambda\in\Part(n)$. Moreover, every Levi subgroup $L$ of $G$ has the form $L_\nu=\GL(\nu_1)\times\GL(\nu_2)\times\cdots$ for some partition $\nu\in\Part(n)$, and it follows that $A_L(x)=1$ for all $x\in\cN_L$. So the condition~\eqref{eqn:big-enough} holds for all $\bk$.

The Weyl group of $\GL(n)$ is the symmetric group $\fS_n$. When $\ell=0$, it is well known that the irreducible $\bk$-representations of $\fS_n$ are parametrized by $\Part(n)$; accordingly, every pair in $\fN_{G,\bk}$ must belong to the principal series, and the generalized Springer correspondence is just the Springer correspondence. In other words, when $\ell=0$, every simple perverse sheaf $\IC(\cO_\lambda,\ubk)$ in $\Perv_G(\cN_G,\bk)$ is a quotient, and hence a direct summand, of the Springer sheaf $\Spr$. It follows that $\GL(n)$ does not have a cuspidal pair unless $n=1$.

Now consider the $\ell>0$ case. The irreducible $\bk$-representations $D^\lambda$ of $\fS_n$ are labelled not by the whole set $\Part(n)$ but by the subset $\Part_\ell(n)$ consisting of partitions that are \emph{$\ell$-regular}, meaning that no part occurs $\ell$ or more times. Under the modular Springer correspondence~\cite{juteau}, $D^\lambda$ corresponds to the pair $(\cO_{\lambda^\tr},\ubk)$ where $\lambda^\tr$ is the transpose partition. Thus the only simple perverse sheaves $\IC(\cO_\lambda,\ubk)$ arising as quotients of $\Spr$ are those where $\lambda$ is \emph{$\ell$-restricted}, i.e.\ $\lambda^\tr$ is $\ell$-regular. (By contrast, using modular reduction from the $\ell=0$ case, one sees that all $\IC(\cO_\lambda,\ubk)$ arise as simple constituents of $\Spr$; see~\cite[Remark~3.2]{genspring1}.) In particular, the principal series is not the only induction series in the modular generalized Springer correspondence unless $\ell>n$.

By Proposition~\ref{prop:distinguished}, the only possible cuspidal pair for $\GL(n)$ is $(\cO_\reg,\ubk)$. In~\cite[Theorem 1]{genspring1} we showed that $(\cO_\reg,\ubk)$ is cuspidal for $\GL(n)$ if and only if $n=\ell^i$ for some nonnegative integer $i$ (our proof there was based on a counting argument, but the claim also follows from Theorem~\ref{thm:sylow} below). As a consequence, the Levi subgroup $L_\nu$ has a cuspidal pair if and only if $\nu$ belongs to the set $\Part(n,\ell)$ of partitions of $n$ into powers of $\ell$, and in this case the unique cuspidal pair is $(\cO_{L_\nu,\reg},\ubk)$. So the set $\fM_{G,\bk}$ is in bijection with $\Part(n,\ell)$, with $\nu\in\Part(n,\ell)$ corresponding to the cuspidal datum $(L_\nu,\cO_{L_\nu,\reg},\ubk)$.

For $\nu\in\Part(n,\ell)$, let $m_{\ell^i}(\nu)$ be the multiplicity of $\ell^i$ as a part of $\nu$; then $n=\sum_{i\geq 0} m_{\ell^i}(\nu)\ell^i$. It is easy to see that
\[
N_G(L_\nu)/L_\nu\cong \prod_{i\geq 0} \fS_{m_{\ell^i}(\nu)} =: \fS_{\sm(\nu)},
\] 
so we have a canonical bijection
\[
\Irr(\bk[N_G(L_\nu)/L_\nu]) \longleftrightarrow \prod_{i\geq 0} \Part_\ell(m_{\ell^i}(\nu)) =: \uPart_{\ell}(\sm(\nu)),
\]
where an element $\boldsymbol \lambda := (\lambda^{(0)},\lambda^{(1)},\lambda^{(2)},\cdots)$ of the right-hand side corresponds to the irreducible representation $D^{\boldsymbol\lambda} := \boxtimes_{i\geq 0} D^{\lambda^{(i)}}$ of $N_G(L_\nu)/L_\nu$. In~\cite[Theorem 3.4]{genspring1} we showed that, under the modular generalized Springer correspondence (see part \eqref{it:mgsc-2} of Theorem~\ref{thm:mgsc}), this irreducible representation corresponds to the pair $(\cO_\lambda,\ubk)$ where
\[
\lambda=\sum_{i\geq 0} \ell^i (\lambda^{(i)})^\tr =: \psico_\nu(\boldsymbol\lambda).
\]

As an example, the pair $(\cO_\reg,\ubk)$ belongs to the induction series associated to the cuspidal datum $(L_{\nu(n)},\cO_{L_{\nu(n)},\reg},\ubk)$, where $\nu(n)$ is the unique element of $\Part(n,\ell)$ for which $m_{\ell^i}(\nu(n))<\ell$ for all $i\geq 0$; in other words, $\nu(n)$ is the partition of $n$ into powers of $\ell$ occurring in the base-$\ell$ expression of $n$. This is a special case of Theorem~\ref{thm:sylow} below.

As might be expected from the proposed theory of modular character sheaves, the above combinatorics is very similar to the combinatorics of induction series for modular representations of finite general linear groups, as in~\cite{dipperdu}.

For concreteness, let us give a complete table for the case $n = 6$, $\ell = 2$.
The first column contains the partitions labelling nilpotent orbits (all local systems are trivial). Each other column represents one induction series. To save space, we write $\GL_k$ instead of $\GL(k)$.
\end{example}

\begin{table}[h]
\[
\begin{array}{|c|c|c|c|c|c|c|}
\hline
L_\nu & \GL_1^6 & \GL_1^4\times\GL_2^1 & \GL_1^2\times\GL_2^2 & \GL_2^3 & \GL_1^2\times\GL_4^1 & \GL_2^1\times\GL_4^1\\
\hline
\fS_{\sm(\nu)} & \fS_6 & \fS_4 \times \fS_1& \fS_2 \times \fS_2 & \fS_3 & \fS_2 \times \fS_1 & \fS_1 \times \fS_1\\
\hline
6 &&&&&& (1,1)\\
51 &&&&& (2,1)&\\
42 &&&& 21&&\\
41^2 && (31,1)&&&&\\
3^2 &&& (2,2)&&&\\
321 & 321&&&&&\\
31^3 && (4,1)&&&&\\
2^3 &&&& 3&&\\
2^21^2 & 42&&&&&\\
21^4 & 51&&&&&\\
1^6 & 6&&&&&\\
\hline
\end{array}
\]
\caption{Modular generalized Springer correspondence for $G = \GL(6)$, $\ell = 2$}
\end{table}

\begin{example} \label{ex:sp2n}
Suppose that $G = \Sp(2n)$ and $\ell = 2$; this case is treated in detail in~\cite{genspring2}.
The $G$-orbits in $\cN_G$ are again classified by Jordan form: we denote them $\cO^\Sp_\lambda$, where $\lambda$ belongs to the set
\[
\Part_{\Sp}(2n)=\{\lambda\in\Part(2n)\,|\,m_{2i+1}(\lambda)\text{ is even for all }i\}.
\]
The distinguished orbits are those for which $\lambda$ consists of distinct even parts \cite[Theorem 8.2.14]{cm}; in other words they are parametrized by $\Part_{2,\Sp}(2n) := \Part_{\Sp}(2n) \cap \Part_2(2n)$.
We will see in Theorem \ref{thm:2cusp} that they all support a cuspidal pair (with the trivial local system).
Hence the Levi subgroups admitting a cuspidal pair are those of the form
\[
L_\nu = \GL(\nu_1) \times \cdots \times \GL(\nu_m) \times \Sp(2n - 2k),
\qquad 0 \le k \le n, \ \nu \in \Part(k,2).
\]

The relative Weyl group $N_{\G(N)}(L_\nu)/(L_\nu)$ is isomorphic to $(\Z/2\Z) \wr \fS_{\sm(\nu)}$. Since $\ell = 2$, all its irreducible representations over $\bk$ factor through the quotient $\fS_{\sm(\nu)}$, and we denote them again by $D^{\boldsymbol\lambda}$, with $\boldsymbol\lambda\in\uPart_2(\sm(\nu))$.

The orbits in $\cN_{L_\nu}$ supporting cuspidal pairs are those of the form
\[
\cO_{[\nu];\mu} := \cO_{(\nu_1)} \times \cdots \times \cO_{(\nu_m)} \times \cO^\Sp_\mu, \qquad
\mu \in \Part_{2,\Sp}(2n - 2k).
\]
Thus, the modular generalized Springer correspondence for $G$ can be regarded as a bijection
\begin{equation*} \label{eqn:sp2n2-genspring2}
\Omega: \bigsqcup_{0 \le k \le n} \  \bigsqcup_{\nu \in \Part(k,2)} \Part_{2,\Sp}(2n-2k)\times\uPart_2(\sm(\nu)) \to \Part_\Sp(2n).
\end{equation*}

 For $\lambda \in \Part(m)$ and $\mu \in \Part(m')$, we define $\lambda\cup\mu\in\Part(m+m')$ to be the partition whose parts are the union of those of $\lambda$ and those of $\mu$; thus, $(\lambda\cup\mu)^\tr=\lambda^\tr+\mu^\tr$.
In \cite[Theorem 9.5]{genspring2} we prove that
\[
\Omega=\bigsqcup_{0 \le k \le n} \  \bigsqcup_{\nu \in \Part(k,2)} \omegaco_{k,\nu},
\]
where $\omegaco_{k,\nu}:\Part_{2,\Sp}(2n-2k)\times\uPart_2(\sm(\nu)) \to \Part_\Sp(2n)$ is defined by
\[
\omegaco_{k,\nu}(\mu,\blambda)=\mu\cup\psico_{k,\nu}(\blambda)\cup\psico_{k,\nu}(\blambda).
\]
Here $\psico_{k,\nu}$ denotes the map $\psico_\nu$ for $\GL(k)$, $\ell=2$.
As an illustration, here is the complete correspondence for $G = \Sp(6)$, $\ell = 2$. Cuspidal objects are denoted by a $*$ symbol.
We abbreviate $\cO_{[\nu];\mu}$
by $\nu_1.\nu_2.\cdots.\nu_m;\mu$.

\begin{table}[h]
\[
\begin{array}{|c|c|c|c|c|c|c|c|}
\hline
L_\nu & \GL_1^3 & \GL_2^1 \times \GL_1^1 & \GL_1^2\times\Sp_2 & \GL_2^1 \times \Sp_2 & \GL_1^1\times\Sp_4 & \multicolumn{2}{c|}{\Sp_6}\\
\hline
\cO_{L_\nu} & 1.1.1 & 2.1 & 1.1;2 & 2;2 & 1 ; 4 & 42 & 6 \\ 
\hline
\fS_{\mathbf m(\nu)} & \fS_3 & \fS_1 \times \fS_1& \fS_2  & \fS_1 & \fS_1 & 1 & 1\\
\hline
6 &&&&&& &*\\
42 &&&&&& *&\\
41^2 &&&&& 1 &&\\
3^2 && (1,1)&&&&&\\
2^3 &&&& 1&&&\\
2^21^2 & 21 &&&&&&\\
21^4 &&& 2 &&&&\\
1^6 & 3 &&&&&&\\
\hline
\end{array}
\]
\caption{Modular generalized Springer correspondence for $G = \Sp(6)$, $\ell = 2$}
\label{table:sp}
\end{table}

\end{example}

\begin{example} \label{ex:g2}
Suppose that $G$ is of type $G_2$; this case is treated in detail in~\cite{genspring3}. There are five $G$-orbits in $\cN_G$, with Bala--Carter labels $G_2$ (the regular orbit), $G_2(a_1)$ (the subregular orbit), $\widetilde{A_1}$, $A_1$ and $\{0\}$. The distinguished orbits are the first two of these. We have $A_G(x)=\fS_3$ for $x\in G_2(a_1)$ and $A_G(x)=1$ for all other $x\in\cN_G$. Since every field is a splitting field for $\fS_3$, and all proper Levi subgroups of $G$ are isomorphic either to $\GL(2)$ or to $\GL(1)^2$, the condition~\eqref{eqn:big-enough} holds for all $\bk$. Since $\fS_3$ has (up to isomorphism) three irreducible $\bk$-representations when $\ell\notin\{2,3\}$ and two when $\ell\in\{2,3\}$, the set $\fN_{G,\bk}$ has seven elements when $\ell\notin\{2,3\}$ and six when $\ell\in\{2,3\}$. Note that $2$ and $3$ are also the primes dividing the order of the Weyl group $W$, which is the dihedral group of order $12$.

When $\ell\notin\{2,3\}$, there are six irreducible $\bk$-representations of $W$, corresponding to six of the seven elements of $\fN_{G,k}$ in the same way as the Springer correspondence (which is the $\ell=0$ case\footnote{Up to the difference in conventions amounting to tensoring with the sign representation; see Footnote~\ref{footnote:sign-twist}.}). The remaining element of $\fN_{G,\bk}$ is cuspidal, namely $(G_2(a_1),\cE_{\mathrm{sign}})$, where $\cE_{\mathrm{sign}}$ is the $G$-equivariant local system on $G_2(a_1)$ corresponding to the sign representation of $\fS_3$.

When $\ell=3$, there are four irreducible $\bk$-representations of $W$, and the modular Springer correspondence between these and four of the six elements of $\fN_{G,\bk}$ was worked out in~\cite[Section 9.1.2]{juteau}. The remaining two elements of $\fN_{G,\bk}$ must be cuspidal, because the proper Levi subgroups of $G$ do not support a cuspidal pair in characteristic $3$, as follows from Example~\ref{ex:gln}. These two cuspidal pairs are $(G_2(a_1),\cE_{\mathrm{sign}})$, in accordance with a general principle that the modular reduction of a characteristic-$0$ cuspidal pair is always cuspidal (see Section~\ref{sec:further})
and $(G_2,\ubk)$, in accordance with Theorem~\ref{thm:sylow} below.

When $\ell=2$, there are only two irreducible $\bk$-representations of $W$, corresponding to the pairs $(\{0\},\ubk)$ and $(\widetilde{A_1},\ubk)$ under the modular Springer correspondence~\cite[Section 9.1.1]{juteau}. There are two cuspidal pairs, as in the $\ell=3$ case: $(G_2(a_1),\ubk)$ and $(G_2,\ubk)$. The other two elements of $\fN_{G,\bk}$ belong to two separate induction series associated to the two $G$-conjugacy classes of Levi subgroups isomorphic to $\GL(2)$, with $(A_1,\ubk)$ attached to the long-root Levi subgroup and $(G_2(a_1),\cE_{\mathrm{refln}})$ attached to the short-root Levi subgroup. (Here $\cE_{\mathrm{refln}}$ is the local system corresponding to the unique $2$-dimensional irreducible $\bk$-representation of $\fS_3$.) For these Levi subgroups $L$, the group $N_G(L)/L$ is $\fS_2$, which indeed has only one irreducible $\bk$-representation.      

We give the tables below. Again, cuspidal objects are denoted by a $*$ (in a single column, to save space). We omit the pair $(\cO_L, \cE_L)$ since it is clear in each case. The notation $\chi_{a,b}$ is the standard notation of the characters of $W$, while $\varphi_{a,b}$ denotes the modular character of $W$ which corresponds to $\chi_{a,b}$ according to the Springer basic set, as explained in  \cite{juteau}. Finally, $\overline \chi_{a,b}$ denotes the modular reduction of $\chi_{a,b}$.
\end{example}

\begin{table}[h]
\[
\begin{array}{ccc}
\begin{array}{|c | c|c|c|c|}
 \multicolumn{5}{c}{\ell = 2} \\[1em]
 \hline
L & T & A_1 & \widetilde A_1 & G_2 \\
\hline
N_G(L)/L & W & \fS_2 & \fS_2 & 1 \\
\hline
G_2 &&&& * \\
G_2(a_1)&&&& * \\
G_2(a_1),\cE_{\mathrm{refln}}&&& 2 &\\
\widetilde A_1 & \varphi_{2,2} = \ov\chi_{2,2} &&&\\
A_1 && 2 &&\\
0 & \varphi_{1,0} = \ov\chi_{1,0} &&& \\
\hline
\end{array}
&&
\begin{array}{|c | c|c|}
\multicolumn{3}{c}{\ell = 3} \\[1em]
\hline
L &    T & G_2
\\
\hline
N_G(L)/L & W & 1\\
\hline
G_2 && *\\
G_2(a_1) &\varphi_{2,1} = \overline \chi_{1,3}''&\\
G_2(a_1),\cE_{\mathrm{sign}} && *\\
\widetilde A_1 & \varphi_{2,2} = \ov\chi_{1,6}&\\
A_1 & \varphi_{1,3}' = \ov\chi_{1,3}' & \\
0 & \varphi_{1,0} = \ov\chi_{1,0}& \\
\hline
\end{array}
\end{array}
\]
\caption{Modular generalized Springer correspondence for $G = G_2$, $\ell = 2, 3$}
\end{table}

\section{Further results}
\label{sec:further}

Theorem~\ref{thm:mgsc} raises the problem of computing the modular generalized Springer correspondence explicitly in terms of the known parametrizations of the sets $\fN_{G,\bk}$ and $\Irr(\bk[N_G(L)/L])$, as was done in the $\ell=0$ case by Lusztig and Spaltenstein~\cite{lusztig,charsh5,lus-spalt,spaltenstein}. Standard principles~\cite[Section 5.3]{genspring2} reduce the problem to the case where $G$ is simply connected and quasi-simple, so one can treat each Lie type in turn; we considered the classical types in~\cite{genspring2} and the exceptional types in~\cite{genspring3}, but in both cases some gaps remain.

It is easy to see that, in a suitable precise sense~\cite[Lemma 3.3(3)]{genspring3}, the modular generalized Springer correspondence depends only on $\ell$, not on the field $\bk$ itself (assuming as always that $\bk$ is big enough to satisfy~\eqref{eqn:big-enough}). Another general fact~\cite[Proposition 7.1]{genspring3}, already seen in Examples~\ref{ex:gln} and~\ref{ex:g2}, is that if $\ell$ does not divide $|W|$, the correspondence is essentially the same as the $\ell=0$ case. The following result illustrates how the correspondence depends on $\ell$ and $|W|$ in general. In this statement, for a Levi subgroup $L$ of $G$, $W_L$ denotes its Weyl group, a parabolic subgroup of $W$.   

\begin{thm} \label{thm:sylow} \cite[Theorem 4.5]{genspring3}
The pair $(\cO_{\reg},\ubk)\in\fN_{G,\bk}$ belongs to the induction series associated to $(L,\cO_{L,\reg},\ubk)$ where $L$ is in the unique $G$-conjugacy class of Levi subgroups that are minimal subject to the condition that $\ell\nmid|W/W_L|$. In particular, $(\cO_{\reg},\ubk)$ belongs to the principal series if and only if $\ell\nmid|W|$, and $(\cO_{\reg},\ubk)$ is cuspidal if and only if no proper parabolic subgroup of $W$ contains an $\ell$-Sylow subgroup of $W$.
\end{thm}
\noindent
This theorem should be compared with the result of Geck--Hiss--Malle~\cite[Theorem 4.2]{ghm} about the modular Steinberg character of a finite reductive group.

The first step towards computing the modular generalized Springer correspondence for a given simply connected quasi-simple group $G$ and prime $\ell$ is the classification of cuspidal pairs, which in the $\ell=0$ case was done by Lusztig~\cite{lusztig}. Our main tools for determining these pairs are the following:
\begin{itemize}
\item
Proposition~\ref{prop:distinguished}, which already eliminates all non-distinguished orbits;
\item
Theorem~\ref{thm:sylow}, which determines when the pair $(\cO_{\reg},\ubk)$ is cuspidal;
\item
\cite[Proposition~2.22]{genspring1}, which says that cuspidality is preserved by \emph{modular reduction}. More precisely, this result implies that if $\K$ is a finite extension of $\Ql$ with ring of integers $\O$ and residue field $\bk$, and if $(\cO,\cE^{\K})$ is a cuspidal pair over $\K$, then for any $\O$-form $\cE^{\O}$ of $\cE^{\K}$ and for any irreducible $G$-equivariant $\bk$-local system $\cF$ on $\cO$ which is a constituent of $\bk \otimes_{\O} \cE^{\O}$, the pair $(\cO,\cF)$ is cuspidal.
\end{itemize}
These tools are sufficient to classify cuspidal pairs in most cases, although some indeterminacies remain for certain bad characteristics in certain exceptional groups; see~\cite[Section 6.4 and Appendix A]{genspring3}. 

In particular, we have the following result about the case where $\ell$ is a good prime for $G$:

\begin{thm}
\label{thm:classification-cuspidals} \cite[Theorem~1.5]{genspring3}
Assume that $G$ is semisimple and simply connected, and that $\ell$ is a good prime for $G$. Then the cuspidal pairs all arise by modular reduction of characteristic-$0$ cuspidal pairs. 
\end{thm}
\noindent
It is important here that we assume that $G$ is semisimple and simply connected. For instance, this statement would fail for $G=\GL(n)$ when $n$ is a power of $\ell$, as we saw in Example~\ref{ex:gln}, and it would therefore also fail for $G=\PGL(n)$ when $n$ is a power of $\ell$. The cuspidal pair $(\cO_\reg,\ubk)$ in these cases really arises by modular reduction of a cuspidal pair $(\cO_\reg,\cE^\K)$ for the group $\SL(n)$, in which $\cE^\K$ has multiplicative order $n$ and thus becomes constant when reduced modulo $\ell$.

By contrast, when $\ell$ is a bad prime for $G$, there can be many more cuspidal pairs than in the $\ell=0$ case, as seen in the following result. 

\begin{thm} \cite[Theorems 7.1 and 8.1]{genspring2} \label{thm:2cusp}
Let $G$ be a simply connected quasi-simple group of type $B$, $C$ or $D$ and let $\ell=2$. Then for \emph{every} distinguished nilpotent orbit $\cO$, the pair $(\cO,\ubk)$, which is the only element of $\fN_{G,\bk}$ supported on $\cO$, is cuspidal.
\end{thm}
\noindent
For groups of these types $B$, $C$ or $D$, this bad characteristic (i.e.\ $\ell=2$) case is in fact the only case where we have explicitly computed the modular generalized Springer correspondence; see~\cite[Theorems 9.5 and 9.7]{genspring2},
and Example \ref{ex:sp2n}.

In~\cite[Section 6]{genspring3} we explained how to determine the number of cuspidal pairs for each simply connected quasi-simple group $G$ of exceptional type and each value of $\ell$; these numbers are given in Table~\ref{tab:exc-cuspidal-pairs}, and again show the differing behaviour of bad primes (recall that the bad primes for a group of exceptional type are $2$, $3$ and, for type $E_8$ only, $5$).

\begin{table}
\[
\begin{array}{|c||c|c|c|c|}
\hline
& \ell=2 & \ell=3 & \ell=5 & \ell\geq 7\\
\hline\hline
E_6 & 4 & 3 & 2 & 2\\
\hline
E_7 & 6 & 3 & 1 & 1\\
\hline
E_8 & 10 & 8 & 5 & 1\\
\hline
F_4 & 4 & 3 & 1 & 1\\
\hline
G_2 & 2 & 2 & 1 & 1\\
\hline
\end{array}
\]
\caption{Number of cuspidal pairs for exceptional simply connected $G$}\label{tab:exc-cuspidal-pairs}
\end{table} 

Theorem~\ref{thm:classification-cuspidals} suggests that one could expect the modular generalized Springer correspondence to have some uniform properties in good characteristic. In fact, a better setting is that of
\emph{rather good} characteristic, i.e.~the case when $\ell$ is good for $G$ and does not divide the order of $Z(G)/Z(G)^\circ$. (See~\cite[\S 2.1]{rather-good} for a discussion of this condition.) If $\ell$ is rather good for $G$ and if $\K$ is a finite extension of $\Ql$ which satisfies condition~\eqref{eqn:big-enough}, then there exists a natural bijection $\fN_{G,\bk} \leftrightarrow \fN_{G,\K}$, so that the following statement is meaningful:

\begin{thm} \label{thm:rather-good} \cite[Theorem~1.1]{rather-good}
If $\ell$ is rather good for $G$, the partition of $\fN_{G,\bk}$ into induction series as in Theorem~{\rm \ref{thm:mgsc}\eqref{it:mgsc-1}} is a refinement of the corresponding partition of $\fN_{G,\K}$, known by the work of Lusztig and Spaltenstein.
\end{thm}

An important property of Lusztig's characteristic-0 generalized Springer correspondence, which plays a crucial role in the theory of character sheaves, is \emph{cleanness}, namely the fact that for every cuspidal pair $(\cO,\cE)$ we have $\IC(\cO,\cE)|_{\overline{\cO} \smallsetminus \cO}=0$. This property fails in general in the modular setting, as can already be seen in the case $G=\mathrm{GL}(2)$, see~\cite[Remark~2.5]{genspring1}. However, Mautner conjectured (in unpublished work) that cleanness holds in rather good characteristic for the cuspidal pairs obtained by modular reduction from a cuspidal pair in characteristic $0$. (By Theorem~\ref{thm:classification-cuspidals}, this covers all cuspidal pairs if $G$ is semisimple and simply connected.)
General principles~\cite[Lemma 2.5]{rather-good} allow us to reduce Mautner's conjecture to the case where $G$ is either a semisimple group of type $A$, or a simply connected quasi-simple group not of type $A$. In~\cite[Theorem~1.3]{rather-good} we have proved the conjecture for groups of type $A$, for groups of exceptional type, for characteristics $\ell$ which do not divide the order of the Weyl group, and also for certain low-rank classical groups. In~\cite[Theorem~1.6]{rather-good} we showed that, provided this conjecture holds, the coarser of the two partitions of $\fN_{G,\bk}$ in Theorem~\ref{thm:rather-good}, namely the one provided by the characteristic-$0$ induction series, induces an orthogonal decomposition of the whole derived category $\Db_G(\cN_G,\bk)$, as Lusztig observed in the $\Qlb$ setting.


\end{document}